\documentclass[10pt]{amsart}
\usepackage{amsmath,amssymb}	
\usepackage{amsthm, amsfonts}
\usepackage{latexsym, mathrsfs}
\usepackage[all,ps,cmtip]{xy}
\usepackage{tabularx}
\usepackage{pict2e}
\usepackage{graphicx}
\usepackage{xypic}
\usepackage{tikz}
\usetikzlibrary{matrix,arrows}
\usepackage{hyperref}
\usepackage[text={145mm, 220mm},centering]{geometry}
\theoremstyle{plain}
\newtheorem{thm}{Theorem}[section]
\newtheorem{lem}[thm]{Lemma}
\newtheorem{prop}[thm]{Proposition}
\newtheorem{cor}[thm]{Corollary}
\newtheorem{claim}[thm]{Claim}
\theoremstyle{definition}
\newtheorem{defn}[thm]{Definition}
\newtheorem{rem}[thm]{Remark}
\newtheorem{example}[thm]{Example}
\newtheorem{quest}[thm]{Question}

\DeclareMathOperator{\CC}{\mathscr{C}}
\DeclareMathOperator{\CN}{\mathbb{C}}

\DeclareMathOperator{\CO}{\mathscr{O}}
\DeclareMathOperator{\CH}{\mathcal{H}}
\DeclareMathOperator{\A}{\mathscr{A}}
\DeclareMathOperator{\CR}{\mathcal{C}}
\DeclareMathOperator{\K}{\mathscr{K}}
\DeclareMathOperator{\BH}{\mathbb{H}}

\allowdisplaybreaks

\numberwithin{equation}{section}

\begin{document}

\title[Blow-up formula of twisted de Rham cohomology]
{On the blow-up formula of twisted de Rham cohomology}

\author[Y. Chen]{Youming Chen}
\address{Y. Chen.: School of Science, Chongqing University of Technology, Chongqing, 400054, People's Republic of China}
\email{youmingchen@cqut.edu.cn}

\author[S. Yang]{Song Yang}
\address{S. Yang.: Center for Applied Mathematics, Tianjin University, Tianjin, 300072, People's Republic of China}
\email{syangmath@tju.edu.cn}

\date{\today}

\begin{abstract}
We derive a blow-up formula for the de Rham cohomology
of a local system of complex vector spaces on a compact complex manifold.
As an application, we obtain the blow-up invariance of $E_{1}$-degeneracy
of the Hodge-de Rham spectral sequence associated with a local system of complex vector spaces.
\end{abstract}



\maketitle

\section{Introduction}

The purpose of this note is to give a blow-up formula for the twisted de Rham cohomology of a local system of complex vector spaces
on a compact complex manifold by the sheaf-theoretical approach.
It is worth to notice that the notion of local system has been extensively studied in complex geometry, especially in Hodge theory.

In this note, we are mainly interested in the twisted de Rham cohomologies of local systems of complex vector spaces on complex manifolds.
Recall that a {\it local system} of complex vector spaces of rank $r$ on a complex manifold $X$
is an abelian sheaf which is locally isomorphic to the constant sheaf of stalk $\CN^{r}$.
Naturally, the constant sheaf $\CN$ is a local system of rank $1$ who admits two important resolutions: the {\it de Rham resolution} and the {\it holomorphic de Rham resolution}.
Then, the famous {\it de Rham theorem} asserts that there exist canonical isomorphisms
\begin{equation}\label{dRthm}
H_{dR}^{l}(X; \CN)
\cong
H^{l}(X, \CN)
\cong
\BH^{l}(X, \Omega_{X}^{\bullet})
\end{equation}
for any $l\in \mathbb{N}$.
Moreover, there is a twisted version for
the de Rham cohomologies of general local systems (cf. Proposition \ref{t-deRham-thm}).

It is well known that the blow-up transformation plays a significant role in complex geometry.
A natural question is to study the blow-up (or more general birational) invariance of the invariants on complex manifolds.
However, the Betti numbers generally are not invariants under the blow up transformations since the blow up centers also have some contributions.
Fortunately, the Betti numbers satisfy the so-called {\it blow-up formula}.
More precisely, suppose $X$ is a compact complex manifold and $Z\subset X$ is a compact complex submanifold of codimension $c\geq 2$, and let  $\widetilde{X}$ be the blow-up of $X$ with the center $Z$.
Then the {\it de Rham blow-up formula} says that there is an isomorphism
\begin{equation}\label{deRhamblowup}
H^{l}(\widetilde{X}, \CN)
\cong
H^{l}(X, \CN)
\oplus \bigoplus_{i=1}^{c-1} H^{l-2i}(Z, \CN).
\end{equation}
It is of importance to notice that many important invariants of complex manifolds hold such a similar formula,
for example,
the Deligne cohomology \cite{BV97},
the Chern classes \cite{GP07}
and the Dolbeault cohomology of holomorphic vector bundles \cite{RYY17,RYY18}, etc.

In the traditional method,
the Thom isomorphism plays an important role in
the proof of the de Rham blow-up formula of complex manifolds (cf. \cite[I, Theorem 7.31]{Voi02});
see also \cite[Appendix A]{RYY17} for another interpretation by using the relative de Rham cohomology.
Based on the twisted version of the second isomorphism in \eqref{dRthm},
we derive the following result without using the Thom isomorphism.

\begin{thm}\label{mainresult}
Let $X$ be a compact complex manifold, and let $\iota: Z\hookrightarrow X$ be a closed complex submanifold of codimension $c\geq 2$.
If $H$ is a local system of complex vector spaces on $X$,
then there is an isomorphism
\begin{equation}\label{mainiso}
H^{l}(\widetilde{X}, \pi^{-1}H)
\cong
H^{l}(X, H)
\oplus \bigoplus_{i=1}^{c-1} H^{l-2i}(Z, \iota^{-1}H)
\end{equation}
where $\pi: \widetilde{X}\rightarrow X$ is the blow-up of $X$ along $Z$.
\end{thm}

For one thing, this result is inspired by the recent work on the blow-up formulae of Dolbeault cohomology and Bott-Chern cohomology \cite{RYY17,RYY18,YY17} and also \cite{ASTT17,Men18b,Ste18} for different discussions on the Dolbeault blow-up formula.
For another, it is a natural non-trivial generalization of the de Rham blow-up formula
since the topological inverse images $ \pi^{-1}\CN\cong\CN$ and $\iota^{-1}\CN\cong\CN$,
and it may also be viewed as a generalization of the blow-up formula for the Morse-Novikov cohomology \cite{YZ15,Men18a}.

Moreover,
for the pair $(X, H)$,
we also have the {\it (twisted) Hodge-de Rham (or Fr\"{o}licher) spectral sequence}
\begin{equation}\label{H-dR-spectral}
E_{1}^{p,q}=H^{q}(X, \Omega_{X}^{p}\otimes_{\CO_{X}} \CH)
\Longrightarrow
H^{p+q}(X, H),
\end{equation}
where $\CH:=H\otimes_{\CN}\CO_{X}$ is the associated locally free sheaf $\CO_{X}$-module with the natural holomorphic flat connection.
If the local system $H$ is the constant sheaf $\CN$ or more general a unitary local system on a compact K\"{a}hler manifold, then the Hodge-de Rham spectral sequence \eqref{H-dR-spectral} degenerates at $E_{1}$.
However, if not,
the $E_{1}$-degeneracy of the Hodge-de Rham spectral sequence generally does not hold even for a compact K\"{a}hler manifold.
This also means that Theorem \ref{mainresult} is non-trivial
even for compact K\"{a}hler manifolds.

As an application of Theorem \ref{mainresult},
we have the following.

\begin{cor}\label{cor1}
Under the hypotheses of Theorem \ref{mainresult}.
The Hodge-de Rham spectral sequence degenerates at $E_{1}$
for $(X, H)$ and $(Z, \iota^{-1}H)$
if and only if
it degenerates at $E_{1}$ for $(\widetilde{X}, \pi^{-1}H)$.
\end{cor}

This result is a natural generalization of \cite[Theorem 1.6]{RYY17}
where the birational invariance of the non-twisted $E_{1}$-degeneracy is also obtained for compact complex threefolds and fourfolds.

The rest of this note is devoted to give a proof of Theorem \ref{mainresult}.
In Sect. \ref{prelim}, we first recall some basic materials on the local systems,
the Gauss-Manin connection and the twisted de Rham theorem,
and then, we review some basics on the blow-up of complex manifolds
and give a description of blow-up of Iwasawa manifold along a smooth curve.
Section \ref{proof} gives the proof of Theorem \ref{mainresult} and Corollary \ref{cor1}.
Finally, in Sect. \ref{finalrem} we give some remarks on the algebraic de Rham cohomology and a further question.

\subsection*{{\bf Notations}}
Assume that $X$ is a compact complex manifold and we fix some notations for later use.
\begin{itemize}
\item $\CN$ the constant sheaf of stalk $\CN$;
\item $\Omega_{X}^{p}$ the sheaf of holomorphic $p$-forms;
\item $\CO_{X}=\Omega_{X}^{0}$ the sheaf of holomorphic functions, i.e., the structure sheaf of complex manifold $X$;
\item $\A_{X}^{l}$ the sheaf of complex-valued smooth $l$-forms;
\item $\A_{X}^{p,q}$ the sheaf of smooth $(p,q)$-forms;
\item $\CC_{X}^{\infty}=\A_{X}^{0}=\A_{X}^{0,0}$ the sheaf of complex-vlaued smooth functions.
\end{itemize}
Suppose $f: Y\rightarrow X$ is a proper holomorphic map between complex manifolds.
We denote by $f_{\ast}(-)$ the direct image of sheaves or currents and $f^{-1}(-)$ the inverse image of sheaves or sets; see for example \cite{Dem12}.


\section{Generalities}\label{prelim}

\subsection{Twisted de Rham cohomology}

In this subsection, we briefly review the twisted de Rham complex of a local system of complex vector spaces which is equivalent to a flat smooth complex vector bundle on a complex manifold.
We mainly refer to \cite[I \S 9.2 and II \S 5.1.1 ]{Voi02}.

Let $X$ be a compact complex manifold of complex dimension $n$ and $H$ be a local system of complex vector spaces of rank $r$.
For $H$, there are two associated locally free sheaves of rank $r$,
$$
\CC_{X}^{\infty}\otimes_{\CN} H
\;\; \textrm{and} \;\;
\CO_{X}\otimes_{\CN}H,
$$
which are $\CC_{X}^{\infty}$-module and $\CO_{X}$-module, respectively.
Since $\CC_{X}^{\infty}$ is an $\CO_{X}$-module, one has
$$
\CC_{X}^{\infty}\otimes_{\CN}H
\cong
(\CO_{X}\otimes_{\CO_{X}}\CC_{X}^{\infty})\otimes_{\CN}H
\cong
\CC_{X}^{\infty}\otimes_{\CO_{X}}(\CO_{X}\otimes_{\CN} H).
$$
Therefore,
without causing confusion,
we may denote by $\CH$ both the locally free sheaves $\CC_{X}^{\infty}\otimes_{\CN}H$ and $\CO_{X}\otimes_{\CN}H$,
and also the associated smooth complex and holomorphic vector bundles as well.

Moreover, in the smooth complex case,
there is a canonical {\it flat smooth connection}, i.e.,
a $\CC_{X}^{\infty}$-linear morphism
\begin{equation}\label{GM-conn-C-infty}
\nabla: \CH\rightarrow\A_{X}^{1}\otimes_{\CC_{X}^{\infty}} \CH
\end{equation}
given by
\begin{equation*}
\nabla\sigma
=
\sum_{i=1}^r d\alpha_{i}\otimes \sigma_{i}=\sum_{i=1}^r (\partial\alpha_{i}\otimes\sigma_{i}+  \bar{\partial}\alpha_{i}\otimes\sigma_{i})
\end{equation*}
for any $\sigma=\sum \alpha_{i}\sigma_{i}\in \CH$ where $\{\sigma_{i}\}$ is a basis of a local trivialization of $H$
and $\{\alpha_{i}\}$ are the corresponding local sections of $\mathscr{C}_{X}^{\infty}$.
Naturally, we can extend the connection \eqref{GM-conn-C-infty} to a $\CC_{X}^{\infty}$-linear morphism
$$
\nabla:\A_{X}^{l}\otimes_{\CC_{X}^{\infty}} \CH \rightarrow \A_{X}^{l+1}\otimes_{\CC_{X}^{\infty}} \CH
$$
by setting $\nabla(\alpha\sigma)=(d\alpha) \sigma+(-1)^{|\alpha|}\alpha\nabla(\sigma)$,
and it decomposes as $(1,0)$-part and $(0,1)$-part
\begin{equation*}
\nabla=\nabla^{1,0}+\nabla^{0,1}.
\end{equation*}
On one hand, the $(0,1)$-part $\nabla^{0,1}$ gives the holomorphic structure,
i.e.,
$(\CC_{X}^{\infty}\otimes_{\CN} H, \nabla^{0,1})=\CO_{X}\otimes_{\CN}H$.
On the other hand, $\nabla^{1,0}$ defines a canonical {\it flat holomorphic connection} on $\CO_{X}\otimes_{\CN}H$ which means a $\CO_{X}$-linear morphism
\begin{equation}\label{hol-connect}
\nabla^{1,0}: \CH\rightarrow \Omega_{X}^{1}\otimes_{\CO_{X}}\CH
\end{equation}
with
\begin{equation*}
\nabla^{1,0}\sigma
=
\sum_{i} \partial\alpha_{i}\otimes\sigma_{i}
\end{equation*}
for any $\sigma=\sum \alpha_{i}\sigma_{i}\in \CH$ where $\{\sigma_{i}\}$ is a basis of a local trivialization of $H$
and $\{\alpha_{i}\}$ are the local holomorphic sections.
Similarly, one can extend the connection \eqref{hol-connect} to a $\CO_{X}$-linear morphism
$\nabla^{1,0}:\Omega_{X}^{l}\otimes_{\CO_{X}} \CH \rightarrow \Omega_{X}^{l+1}\otimes_{\CO_{X}} \CH$ with the Leibniz rule.

In both cases, the canonical flat connections \eqref{GM-conn-C-infty} and \eqref{hol-connect} are called {\it Gauss-Manin connections}.
In summary,
the smooth complex or holomorphic vector bundle associated with $H$ is naturally equipped with a canonical flat connection,
and actually, the converse also holds.
More exactly,
we have the following well-known fact (cf. \cite[I, Proposition 9.11]{Voi02}).

\begin{lem}[Deligne]
The map $H\mapsto (\CH, \nabla)$ (resp. $H\mapsto (\CH, \nabla^{1,0})$) gives a one-to-one correspondence between
the isomorphism classes of local systems of complex vector spaces
and the isomorphism classes of smooth complex (resp. holomorphic) vector bundles endowed with a flat connection.
\end{lem}

Moreover, the Gauss-Manin connections \eqref{GM-conn-C-infty} and \eqref{hol-connect} may naturally be extended to the differentials, denoted by
$$
d_{\nabla}:
\A_{X}^{l}\otimes_{\CC_{X}^{\infty}} \CH
\longrightarrow
\A_{X}^{l+1}\otimes_{\CC_{X}^{\infty}}  \CH
$$
and
$$
\partial_{\nabla}:
\A_{X}^{p,q}\otimes_{\CC_{X}^{\infty}}\CH
\longrightarrow
\A_{X}^{p+1,q}\otimes_{\CC_{X}^{\infty}}  \CH,
$$
with the Leibniz rule.
For simplicity, from now on we drop the subscript $\nabla$ and we denote by $\A_{X}^{l}(\CH):=\A_{X}^{l}\otimes_{\CC_{X}^{\infty}} \CH$ and $\A_{X}^{p,q}(\CH):=\A_{X}^{p,q}\otimes_{\CC_{X}^{\infty}} \CH$.
Then we have differential
$$
d=\partial+\bar{\partial}: \A_{X}^{l}(\CH)\rightarrow \A_{X}^{l+1}(\CH).
$$
Consequently, given a local system $H$ of complex vector spaces,
there is a natural double complex
$$
\big(\A_{X}^{\bullet,\bullet}(\CH), \partial, \bar{\partial} \big)
$$
with the total complex being the {\it twisted de Rham complex}
$$
\xymatrix{
0\ar[r]&\A_{X}^{0}(\CH)\ar[r]^{d}& \A_{X}^{1}(\CH)\ar[r]^{d}& \A_{X}^{2}(\CH)
\ar[r]^{\;\;\; d}&\cdots \ar[r]^{d\;\;}&\A_{X}^{2n}(\CH) \ar[r]& 0,
}
$$
which is a fine resolution of $H$.

\begin{defn}[Twisted de Rham cohomology]
The {\it twisted de Rham cohomology of the local system $H$} is defined as
$$
H_{dR}^{l}(X; H)
:=
\frac{\ker\big(A^{l}(X; H)\stackrel{d}\longrightarrow A^{l+1}(X; H)\big)}{\mathrm{im} \big(A^{l-1}(X; H)\stackrel{d}\longrightarrow A^{l}(X; H)\big)}
$$
for any $0\leq l\leq 2n$, where $A^{l}(X; H):=\Gamma(X, \A_{X}^{l}(\CH))$ is the space
of $\CH$-valued differential $l$-forms.
\end{defn}

Similarly, we denote by $\Omega_{X}^{p}(\CH):=\Omega_{X}^{p}\otimes_{\CO_{X}} \CH$.
Then there is a {\it twisted holomorphic de Rham complex}
$$
\xymatrix{
0\ar[r]& \CO_{X}(\CH) \ar[r]^{\partial}& \Omega_{X}^{1}(\CH) \ar[r]^{\partial}&\cdots \ar[r]^{\partial}& \Omega_{X}^{n}(\CH)\ar[r]& 0.
}
$$
which is also a resolution of $H$.
As a result, one has the twisted de Rham theorem.

\begin{prop}[Twisted de Rham theorem]\label{t-deRham-thm}
$$
H_{dR}^{l}(X; H)
\cong
H^{l}(X, H)
\cong
\BH^{l}(X, \Omega_{X}^{\bullet}(\CH)).
$$
\end{prop}

Finally, we give a current description of twisted de Rham cohomology for the pair $(X, H)$.
We refer the reader to \cite{Dem12} for the current theory.
We denote $\CR_{X}^{l}(\CH)$ the sheaf of $\CH$-valued currents of degree $l$.
Then the de Rham cohomology of $\CH$-valued currents
$$
H_{\CR}^{l}(X; H)
:=
\frac{\ker\big(\CR^{l}(X; H)\stackrel{d}\longrightarrow \CR^{l+1}(X; H)\big)}{\mathrm{im} \big(\CR^{l-1}(X; H)\stackrel{d}\longrightarrow \CR^{l}(X; H)\big)}
$$
for any $0\leq k\leq n$, where $\CR^{l}(X; H):=\Gamma(X, \CR_{X}^{l}(\CH))$ is the space of $\CH$-valued currents of degree $l$ with the natural differential induced by $d$.
From definition, there is a natural inclusion
$j:A^{\bullet}(X; H)\hookrightarrow\CR^{\bullet}(X; H)$
which induces an isomorphism
\begin{equation}\label{currentform}
j:H_{dR}^{l}(X; H)
\stackrel{\cong}\longrightarrow
H_{\CR}^{l}(X; H).
\end{equation}
Moreover, $\CR_{X}^{\bullet}(\CH)$ is also a fine resolution of $H$.
Hence there is a diagram
$$
\xymatrix@C=0.5cm{
& H^{l}(X, H)  \ar[ld]_{\simeq} \ar[rd]^{\simeq} & \\
H_{dR}^{l}(X; H)\ar[rr]_{\simeq}&& H_{\CR}^{l}(X; H)
}
$$

\subsection{Blow-up}

In this subsection, we collect some basic facts on the blow-ups of complex manifolds along its closed complex submanifolds.
We refer the interested readers to \cite[I, \S 3.3.3]{Voi02} for more detailed discussions.

Let $X$ be a complex manifold of complex dimension $n\geq 2$,
and let $Z\subset X$ be a closed complex submanifold of complex dimension $m=n-c$ ($c\geq 2$).
Theoretically, there is a unique complex manifold $\mathrm{Bl}_{Z}X$ and a holomorphic map $\pi:\mathrm{Bl}_{Z}X \longrightarrow X$,
which is called the {\it blow-up of $X$ along $Z$}, such that
\begin{enumerate}
\item  $\pi$ is proper and the restriction $\pi: \mathrm{Bl}_{Z}X \setminus \pi^{-1}(Z)\longrightarrow X\setminus Z$ is isomorphic;
\item  $E:=\pi^{-1}(Z)$ is a hypersurface, which is called the \textit{exceptional divisor},
            of $\mathrm{Bl}_{Z}X$.
           Furthermore, $E$ is isomorphic to the projective bundle $\mathbb{P}(\mathcal{N}_{Z/X})$
of the normal bundle $\mathcal{N}_{Z/X}$ of $Z$ in $X$.
\end{enumerate}
Moreover, if $X$ is compact, so is $\mathrm{Bl}_{Z}X$.

For convenience, we sketch the local construction of blow-up
and then one can glue them into the global blow-up (cf. \cite[I, Lemma 3.22]{Voi02}).
In fact, one defines the blow-up $\mathrm{Bl}_{\mathbb{C}^{m}} \mathbb{C}^{n}$ of $\mathbb{C}^{n}$ along
its closed complex submanifold
$\mathbb{C}^{m}:=\{z_{1}=z_{2}=\cdots=z_{n-m}=0\}$ as follows:
$$
\mathrm{Bl}_{\mathbb{C}^{m}} \mathbb{C}^{n}:=
\{ (w, z) \in \mathbb{P}^{n-m-1}\times \mathbb{C}^{n}\mid w_{i}z_{j}=z_{i}w_{j}, \forall i, j\leq n-m \rangle \}.
$$
It is a complex manifold
and the projection $\pi: \mathrm{Bl}_{\mathbb{C}^{m}} \mathbb{C}^{n}\longrightarrow \mathbb{C}^{n}$ is holomorphic,
$E= \pi^{-1}(\mathbb{C}^{m})\cong \mathbb{C}^{m}\times \mathbb{P}^{n-m-1}$.
Then, the local blow-up diagram is as follows:
\begin{equation*}
\xymatrix{
 \mathbb{C}^{m}\times \mathbb{P}^{n-m-1} \ar[d]_{\rho} \ar@{^{(}->}[r]^{\;\; \tilde{\iota}} & \mathrm{Bl}_{\mathbb{C}^{m}} \mathbb{C}^{n}\ar[d]^{\pi}\\
\mathbb{C}^{m} \ar@{^{(}->}[r]^{\iota} & \mathbb{C}^{n}.
}
\end{equation*}

Now, we consider a concrete example.

\begin{example}[Iwasawa manifold]\label{ex-Iwasawa}
We consider the blow-up of Iwasawa manifold $\mathbf{I(3)}$ along a smooth curve $C$.
Denote by $\mathbf{H}(3; \mathbb{C})$ the Heisenberg Lie group
$
\Big\{ \small{\left(\begin{array}{ccc}
1 & z_{1} & z_{3}\\
0 & 1  & z_{2} \\
0 & 0 & 1
\end{array}
\right)}\mid z_{i}\in \mathbb{C} \Big\}\subset \mathbf{Gl}(3; \mathbb{C}).
$
As complex manifolds, $\mathbf{H}(3; \mathbb{C})$ is isomorphic to $\mathbb{C}^{3}$.
Consider the discrete group
$\mathbf{G(3)}:=\mathbf{Gl}(3; \mathbb{Z}[i])\cap \mathbf{H}(3; \mathbb{C})$,
where $\mathbb{Z}[i]=\{a+bi\mid a,b\in \mathbb{Z}\}$ is the Gaussian integers.
Then the left multiplication gives a natural $\mathbf{G(3)}$-action on $\mathbf{H}(3; \mathbb{C})$.
Correspondingly, there is a faithful $\mathbf{G(3)}$-action on $\mathbb{C}^{3}$ given by
$$
(g_{1}, g_{2}, g_{3})\cdot(z_{1}, z_{2}, z_{3}):=(z_{1}+g_{1}, z_{2}+g_{2}, z_{3}+g_{1}z_{2}+g_{3}),
$$
where $g_{1}, g_{2}, g_{3}\in \mathbb{Z}[i]$.
This $\mathbf{G(3)}$-action defines a monomorphism $\varphi: \mathbf{G(3)}\to \mathbf{Aff}(\mathbb{C}^{3})$
where $\mathbf{Aff}(\mathbb{C}^{3})$ is the affine transformation group of $\mathbb{C}^{3}$.
Therefore,
this $\mathbf{G(3)}$-action is properly discontinuous.
Consequently, the associated $\mathbf{G(3)}$-quotient space
$$
\mathbf{I(3)}:=\mathbb{C}^{3}/ \mathbf{G(3)}
$$
is a compact complex threefold, which is called {\it Iwasawa manifold}.

Moreover, there is a natural holomorphic map $p: \mathbf{I(3)} \longrightarrow \mathbf{T}_{\mathbb{C}}^{2}$
which is induced by the projection $(z_{1}, z_{2}, z_{3})\longmapsto (z_{1}, z_{2})$ where
$
\mathbf{T}_{\mathbb{C}}^{2}:=\mathbb{C}/ \mathbb{Z}[i] \times\mathbb{C}/ \mathbb{Z}[i].
$
The fibers of $p$ are all isomorphic to the complex torus $\mathbf{T}_{\mathbb{C}}=\mathbb{C}/ \mathbb{Z}[i]$.
We denote $C:=p^{-1}([0, 0])=p^{-1}([c_{1}, c_{2}])$ where $c_{1}, c_{2}\in \mathbb{Z}[i]$.
It is a  smooth curve in $\mathbf{I(3)}$.
We claim that $C$ is also a $\mathbf{G(3)}$-quotient space.
To this end, we take
$$
\Gamma: =\mathbb{Z}[i]\times \mathbb{Z}[i] \times \mathbb{C}=\{(c_{1}, c_{2}, z_{3}) \mid, c_{1}, c_{2}\in \mathbb{Z}[i],  z_{3}\in \mathbb{C}\}
$$
as a closed complex submanifold of $\mathbb{C}^{3}$.
It is a direct check that $\Gamma$ is $\mathbf{G(3)}$-equivariant.
Therefore, by the definition,
$
C= \Gamma /\mathbf{G(3)}.
$

Now we consider the blow-up of $\mathbb{C}^{3}$ along $\Gamma$.
One has the following diagram
\begin{equation*}
\xymatrix{
\Gamma \times \mathbb{P}^{1} \ar[d]_{\rho} \ar@{^{(}->}[r]^{\;\; \tilde{\iota}} & \mathrm{Bl}_{\Gamma} \mathbb{C}^{3}\ar[d]^{\pi}\\
\Gamma \ar@{^{(}->}[r]^{\iota} & \mathbb{C}^{3},
}
\end{equation*}
where $\mathrm{Bl}_{\Gamma} \mathbb{C}^{3}\cong \mathbb{Z}[i]\times \mathbb{Z}[i] \times \mathrm{Bl}_{\mathbb{C}} \mathbb{C}^{3}$ as complex manifolds since $\mathbb{Z}[i]$ is discrete in $\mathbb{C}$.
It not difficult to see that the $\mathbf{G(3)}$-action can be lifted into the blow-up diagram.
Since the blow-up complex manifold is unique,
it yields
$$
\mathrm{Bl}_{C}\mathbf{I(3)}\cong \mathrm{Bl}_{\Gamma} \mathbb{C}^{3}/ \mathbf{G(3)}.
$$
\end{example}

\section{Proof}\label{proof}

Let $X$ be a compact complex manifold of dimension $n\geq 2$
and $\iota: Z\hookrightarrow X$ be its closed complex submanifold
of codimension $c\geq 2$.
Let $\pi:\tilde{X}\rightarrow X$ be the blow-up of $X$ along the center $Z$.
Then, we have the blow-up diagram
\begin{equation}\label{blowup-diag}
\xymatrix{
E \ar[d]_{\rho} \ar@{^{(}->}[r]^{\tilde{\iota}} & \widetilde{X}\ar[d]^{\pi}\\
 Z \ar@{^{(}->}[r]^{\iota} & X.
}
\end{equation}
where the exceptional divisor $\rho: E \rightarrow Z$ is the projective bundle of the normal bundle of $Z$ in $X$, and $\tilde{\iota}$ is closed embedding of $E$ in $\tilde{X}$.

From now on, let $H$ be a local system of complex vector spaces on $X$.
Then we have
$$
\widetilde{\CH}:=
\pi^{\ast}\CH
=
\pi^{-1}(H\otimes_{\CN} \CO_{X})\otimes_{\pi^{-1}\CO_{X}}\CO_{\widetilde{X}}
\cong
(\pi^{-1}H\otimes_{\pi^{-1}\CN} \pi^{-1}\CO_{X})\otimes_{\pi^{-1}\CO_{X}}\CO_{\widetilde{X}}
\cong
\pi^{-1}H\otimes_{\CN}\CO_{\widetilde{X}}.
$$
This means that the inverse images of the locally free sheaves and the local systems are compatible.

\subsection{Proof of Theorem \ref{mainresult}}
We use a technical notion of relative Dolbeault sheaf
which is introduced in \cite{RYY17,RYY18,YY17}.
Recall that the pullback of holomorphic $p$-forms determines
a natural surjective sheaf morphism
$\iota^{\ast}:\Omega_{X}^{p}\rightarrow\iota_{\ast}\Omega_{Z}^{p}$
of $\CO_{X}$-modules on $X$.
Therefore, there is a short exact sequence of $\CO_{X}$-modules
\begin{equation}\label{rel-Dol-sequ}
\xymatrix@=0.5cm{
  0 \ar[r] & \K_{X,Z}^{p} \ar[r]^{} & \Omega_{X}^{p} \ar[r]^{\iota^{\ast}} &  \iota_{\ast}\Omega_{Z}^{p} \ar[r] & 0,}
\end{equation}
and we call the kernel sheaf of $\iota^{\ast}$,
denoted by $\K_{X,Z}^{p}$, the {\it $p$-th relative Dolbeault sheaf} of $X$ with respect to $Z$.
Tensoring \eqref{rel-Dol-sequ} with $\CH$,
we obtain a new short exact sequence of complexes of sheaves
\begin{equation}\label{rel-Dol-complex1}
\xymatrix@=0.5cm{
0\ar[r]& \K_{X,Z}^{\bullet}(\CH)
\ar[r]& \Omega_{X}^{\bullet}(\CH)
\ar[r]& \iota_{\ast}\Omega_{Z}^{\bullet}(\iota^{\ast}\CH)
\ar[r]& 0
}
\end{equation}
since the pullback operator commutes with the operator $\partial$.

Similarly, for the pair $(\widetilde{X}, E)$,
we also have a short exact sequence
\begin{equation}\label{rel-Dol-complex2}
\xymatrix@=0.5cm{
0\ar[r]& \K_{\widetilde{X},E}^{\bullet}(\widetilde{\CH})
\ar[r]& \Omega_{\widetilde{X}}^{\bullet}(\widetilde{\CH})
\ar[r]& \tilde{\iota}_{\ast}\Omega_{E}^{\bullet}(\tilde{\iota}^{\ast}\widetilde{\CH})
\ar[r]& 0.
}
\end{equation}

Consider the long exact sequences of hypercohomology of the above two short exact sequences \eqref{rel-Dol-complex1} and \eqref{rel-Dol-complex2}.
Since the direct image $\iota_{\ast}$ and $\tilde{\iota}_{\ast}$ are exact functors,
one has that
 $$
 \BH^{l}(Z, \Omega_{Z}^{\bullet}(\iota^{\ast}\CH))
 \cong
 \BH^{l}(X, \iota_{\ast}\Omega_{Z}^{\bullet}(\iota^{\ast}\CH))
 $$
 and
 $$
 \BH^{l}(E,\Omega_{E}^{\bullet}(\tilde{\iota}^{\ast}\widetilde{\CH}))
 \cong
 \BH^{l}(\widetilde{X},\tilde{\iota}_{\ast}\Omega_{E}^{\bullet}(\tilde{\iota}^{\ast}\widetilde{\CH})).
 $$
Then the blow-up diagram \eqref{blowup-diag} and the pullback of differential forms naturally induce a commutative ladder of long exact sequences of hypercohomologies
\begin{equation*}
\xymatrix@C=0.5cm{
\ar[r]^{} & \BH^{l}(X,\K_{X,Z}^{\bullet}(\CH)) \ar[d]_{\pi^{\ast}} \ar[r]^{} &\BH^{l}(X,\Omega_{X}^{\bullet}(\CH)) \ar[d]_{\pi^{\ast}} \ar[r]^{} & \BH^{l}(Z, \Omega_{Z}^{\bullet}(\iota^{\ast}\CH))\ar[d]_{\rho^{\ast}} \ar[r]^{} & \BH^{l+1}(X,\K_{X,Z}^{\bullet}(\CH))\ar[d]_{\pi^{\ast}} \ar[r]^{} &  \\
\ar[r] & \BH^{l}(\widetilde{X},\K_{\widetilde{X},E}^{\bullet}(\widetilde{\CH})) \ar[r]^{} &
  \BH^{l}(\widetilde{X},\Omega_{\widetilde{X}}^{\bullet}(\widetilde{\CH})) \ar[r]^{} &
   \BH^{l}(E,\Omega_{E}^{\bullet}(\tilde{\iota}^{\ast}\widetilde{\CH}))\ar[r]^{} &
   \BH^{l+1}(\widetilde{X},\K_{\widetilde{X},E}^{\bullet}(\widetilde{\CH})) \ar[r] & }
\end{equation*}

Moreover, by the twisted de Rham theorem (see Proposition \ref{t-deRham-thm}) and the blow-up commutative diagram \eqref{blowup-diag},
the above commutative ladder becomes
\begin{equation}\label{main-commut-ladder}
\xymatrix@C=0.5cm{
 \ar[r]^{} & \BH^{l}(X,\K_{X,Z}^{\bullet}(\CH)) \ar[d]_{\pi^{\ast}} \ar[r]^{} & H^{l}(X, H) \ar[d]_{\pi^{\ast}} \ar[r]^{} & H^{l}(Z, \iota^{-1}H)\ar[d]_{\rho^{\ast}} \ar[r]^{} & \BH^{l+1}(X,\K_{X,Z}^{\bullet}(\CH))\ar[d]_{\pi^{\ast}} \ar[r]^{} &  \\
\ar[r] & \BH^{l}(\widetilde{X},\K_{\widetilde{X},E}^{\bullet}(\widetilde{\CH})) \ar[r]^{} &
  H^{l}(\widetilde{X}, \pi^{-1}H) \ar[r]^{} &
  H^{l}(E, \rho^{-1}(\iota^{-1}H))\ar[r]^{} &
   \BH^{l+1}(\widetilde{X},\K_{\widetilde{X},E}^{\bullet}(\widetilde{\CH})) \ar[r] & }
\end{equation}

To finish the proof, we need the following three lemmas which will be proved later.
The first two lemmas are known for experts,
and we believe that the third one is new.

\begin{lem}\label{injective-lem}
The morphism
$$
\pi^{\ast}:
H^{l}(X, H)
\longrightarrow
H^{l}(\widetilde{X}, \pi^{-1}H)
$$
is injective.
\end{lem}

\begin{lem}\label{bundle-formula}
There is a canonical isomorphism
$$
\bigoplus_{i=0}^{c-1} H^{l-2i}(Z, \iota^{-1}H)
\xrightarrow{\sum\limits_{i=0}^{c-1} h^{i} \wedge \rho^{\ast}(-)}
H^{l}(E, \rho^{-1}(\iota^{-1}H)),
$$
where $h=c_{1}(\CO_{E}(1))$ is the first Chern class of the relative tautological line bundle.
\end{lem}

\begin{lem}\label{key-lem}
The morphism
\begin{equation}\label{key-lem-iso}
\pi^{\ast}:
\BH^{l}(X,  \K_{X,Z}^{\bullet}(\CH))
\stackrel{\simeq}\longrightarrow
\BH^{l}(\widetilde{X},  \K_{\widetilde{X},E}^{\bullet}(\widetilde{\CH}))
\end{equation}
is isomorphic.
\end{lem}

Temporarily admitting the above three lemmas,
we may complete the proof of Theorem \ref{mainresult}.
In the diagram \eqref{main-commut-ladder},
following Lemma \ref{key-lem} the first and the fourth column maps are isomorphisms
and following Lemmas \ref{injective-lem} and \ref{bundle-formula} the second and third are injective.
Then applying the diagram-chasing in \eqref{main-commut-ladder},
we obtain that the cokernel of
$H^{l}(X, H)
\stackrel{\pi^{\ast}}\longrightarrow
H^{l}(\widetilde{X}, \pi^{-1}H)$ is canonical isomorphic to the cokernel
 $H^{l}(Z, \iota^{-1}H)
 \stackrel{\rho^{\ast}} \longrightarrow
 H^{l}(E, \rho^{-1}(\iota^{-1}H))$.
Consequently, by Lemma \ref{bundle-formula} again
we obtain
\begin{eqnarray*}
H^{l}(\widetilde{X}, \pi^{-1}H)
&\cong&
H^{l}(X, H)
 \oplus H^{l}(E, \rho^{-1}(\iota^{-1}H)) /H^{l}(Z, \iota^{-1}H) \\
&\cong&
H^{l}(X, H)
 \oplus \bigoplus_{i=1}^{c-1} H^{l-2i}(Z, \iota^{-1}H).
\end{eqnarray*}
This completes the proof of Theorem \ref{mainresult}.

\subsection{Proof of Lemma \ref{injective-lem}}
Likewise to \eqref{currentform},
for the pair $(\widetilde{X},  \pi^{-1}H)$ one has a canonical isomorphism
\begin{equation*}
\tilde{j}:H_{dR}^{l}(\widetilde{X}; \pi^{-1}H)
\stackrel{\simeq}\longrightarrow
H_{\CR}^{l}(\widetilde{X};  \pi^{-1}H).
\end{equation*}
Moreover, since $\pi$ is a proper map, one has the direct image of the currents
$$
\pi_{\ast}:
H_{\CR}^{l}(\widetilde{X};  \pi^{-1}H)
\longrightarrow
H_{\CR}^{l}(X;  H).
$$
Then we have the following diagram
\begin{equation*}
\xymatrix@C=0.5cm{
 H_{dR}^{l}(X; H)  \ar[d]_{\pi^{\ast}} \ar[r]^{j}_{\simeq} & H^{l}_{\CR}(X; H)\\
 H_{dR}^{l}(\widetilde{X};  \pi^{-1}H)  \ar[r]^{\tilde{j}}_{\simeq} &
 H_{\CR}^{l}(\widetilde{X};  \pi^{-1}H)  \ar[u]_{\pi_{\ast}} }
\end{equation*}
Moreover the blow-up morphism $\pi$ has finite generic fiber of cardinal $1$,
i.e., of degree $1$,
one can show that
$$
\pi_{\ast}\tilde{j}\pi^{\ast}=j
$$
by following the totally same step as \cite[Lemma 2.3]{Wel74} or \cite[Theorem 12.9]{Dem12}.
The fact is essentially due to the condition that $\pi$ is biholomorphic outside sets of Lebesgue measure zero.

Therefore if $\pi^{\ast}(\alpha)=0$, we have that $j(\alpha)=\pi_{\ast}\tilde{j}\pi^{\ast}(\alpha)=0$.
However $j$ is an isomorphism, thus $\alpha=0$ and this means $\pi^{\ast}$ is injective.

\subsection{Proof of Lemma \ref{bundle-formula}}
This is in fact a twisted version of the Leray-Hirsch lemma, and it can be proved by using the Mayer-Vietoris sequence as \cite[Theorem 5.11]{BT82}.
Here, we give a different interpretation following an idea of Deligne \cite{Del68}.

Let $h:=c_{1}(\CO_{E}(1))$ be the first Chern class of the relative tautological line bundle.
and $h^{i}:=\wedge^{i}h \in H^{2i}(E, \CN)$,
$i=0, 1,\ldots,c-1$.
Then the restriction of those classes to each fiber $E_{z}\cong \mathbb{P}^{c-1}$ forms a basis of
$H^{\bullet}(E_{z} ,\CN)$.
Given a local system $L$ of complex vector spaces of rank $r$ on $Z$,
there is a natural morphism of sheaf complexes,
$$
\A_{Z}^{\bullet}(\mathcal{L})[-2i]
\xrightarrow{h^{i}\wedge \rho^{\ast}(-)}
\rho_{\ast}\A_{E}^{\bullet}(\rho^{\ast}\mathcal{L}),
$$
for each $h^{i}$.
Taking direct sums,
we obtain a morphism of complexes of sheaves
\begin{equation}
\Phi:
\bigoplus_{i=0}^{c-1} \A_{Z}^{\bullet}(\mathcal{L})[-2i]
\xrightarrow{\sum_{i=0}^{c-1} h_{i} \wedge \rho^{\ast}(-)}
\rho_{\ast}\A_{E}^{\bullet}(\rho^{\ast}\mathcal{L}).
\end{equation}
Moreover, $\Phi$ is defined up to homotopy by the cohomology class $[h^{i}]$.
In fact, if $h^{i}-h'=dv$ then one has
$
h^{i}\wedge \rho^{\ast}-h'\wedge \rho^{\ast}
=dv\wedge \rho^{\ast}=d(v\wedge \rho^{\ast})-v\wedge d\rho^{\ast}.
$
This means $h^{i}\wedge \rho^{\ast}$ and $h'\wedge \rho^{\ast}$ are homotopic equivalent.

Moreover $\Phi$ is a quasi-isomorphism.
To prove this fact, we consider the induced morphism of cohomology sheaves
$$
H^{t}(\Phi):
\bigoplus_{i=0}^{c-1} H^{t-2i}(\A_{Z}^{\bullet}(\mathcal{L}))
\xrightarrow{\sum\limits_{i=0}^{c-1} h^{i} \wedge \rho^{\ast}(-)}
H^{t}(\rho_{\ast}\A_{E}^{\bullet}(\rho^{\ast}\mathcal{L})).
$$
Since $\rho$ is proper and $\A_{E}^{\bullet}(\rho^{\ast}\mathcal{L})$ is a fine resolution of $\rho^{-1}L$,
on the stalk one has
$$
\big(H^{t}(\rho_{\ast}\A_{E}^{\bullet}(\rho^{\ast}\mathcal{L})\big)_{z}
=
H^{t}(E_{z}, \rho^{-1}L|_{E_{z}})
=H^{t}(E_{z}, \CN)^{\oplus r}
\cong
\begin{cases}
0,  \;\;\; t\  \textrm{is odd}; \\
\CN^{r}, \; t \  \textrm{is  even}.
\end{cases}
$$
Similarly, $\A_{Z}^{\bullet}(\mathcal{L})$ is a fine resolution of $L$;
therefore, it suffices to prove the isomorphism on stalks
$$
H^{2i}(\Phi)_{z}:
L_{z}
\xrightarrow{h^{i}|_{E_{z}} \wedge \rho^{\ast}(-)}
H^{2i}(\rho_{\ast}\A_{E}^{\bullet}(\rho^{\ast}\mathcal{L}))_{z}.
$$
It is indeed an isomorphism since $L_{z}=\CN^{r}$ and $H^{2i}(\Phi)_{z}(1)=h^{i}|_{E_{z}}$ is a basis of $H^{2i}(E_{z}, \CN)$.

\begin{rem}
Using the same idea, another interpretation of the (bundle-valued) Dolbeault projective bundle formula \cite{RYY17,RYY18} may be obtained.
In fact, for any holomorphic vector bundle $\mathcal{L}$ on $Z$,
there is a similar morphism of complexes of sheaves
$$
\bigoplus_{i=0}^{c-1} \A_{Z}^{p-i,\bullet}(\mathcal{L})[-i]\xrightarrow{\sum\limits_{i=0}^{c-1}h_{i}\wedge \rho^{\ast}(-)}
\rho_{\ast}\A_{E}^{p,\bullet}(\rho^{\ast}\mathcal{L}),
$$
which is indeed a quasi-isomorphism.
As a matter of fact the quasi-isomorphism is due to the canonical isomorphism
$$
\rho^{\ast}:
\Omega_{Z}^{p-i}(\mathcal{L})
\stackrel{\simeq}\longrightarrow
R^{i}\rho_{\ast}\Omega_{E}^{p}(\rho^{\ast}\mathcal{L}),
$$
for $0\leq i\leq c-1$.
\end{rem}


\subsection{Proof of Lemma \ref{key-lem}}
Recall
$
\K_{X,Z}^{p,q}(\CH)
:=
\ker\big(\A_{X}^{p,q}(\CH) \stackrel{\iota^{\ast}}\longrightarrow \iota_{\ast}\A_{Z}^{p,q}(\iota^{\ast}\CH)\big)
$
is the relative $(p,q)$-Dolbeault sheaf.
Correspondingly, $\K_{X,Z}^{p,\bullet}(\CH)$ is a fine resolution of $\K_{X,Z}^{p}(\CH)$.

Let $K^{p,q}(\CH)=\Gamma(X, \K_{X,Z}^{p,q}(\CH))$.
The hypercohomology of $\K_{X,Z}^{\bullet}(\CH)$ may be viewed
in terms of a double complex $K^{\bullet,\bullet}(\CH)$ with the differentials
$$
\partial:
K^{p,q}(\CH)\rightarrow  K^{p+1,q}(\CH),\;\;\;\;\;\;
\bar{\partial}:
K^{p,q}(\CH)\rightarrow  K^{p,q+1}(\CH).
$$
Denote the total complex by
$$
K^{l}(\CH):=\bigoplus_{p+q=l}K^{p,q}(\CH).
$$
Then, there is a canonical isomorphism
$
\BH^{l}(X,\K_{X,Z}^{\bullet}(\CH))\cong
H^{l}(K^{\bullet}(\CH)).
$
Similarly, there is a double complex $K^{\bullet,\bullet}(\widetilde{\CH})$ associated with the pair $(\widetilde{X},E)$ and a canonical isomorphism
$
\mathbb{H}^{l}(\widetilde{X},\K_{\widetilde{X},E}^{\bullet}(\widetilde{\CH}))\cong
H^{l}(K^{\bullet}(\widetilde{\CH})).
$
Hence, the morphism \eqref{key-lem-iso} can be rewritten as the morphism of the total cohomology
\begin{equation}\label{morph-total}
\pi^{\ast}:
H^{l}(K^{\bullet}(\CH))
\rightarrow
H^{l}(K^{\bullet}(\widetilde{\CH})).
\end{equation}
Therefore, to conclude the proof it suffices to show that the morphism \ref{morph-total} is isomorphic.
Consider the following two bounded double complexes
$$
(K^{\bullet,\bullet}(\CH); \bar{\partial}, \partial)
\;\;
\textrm{and}\;\;
(K^{\bullet,\bullet}(\widetilde{\CH}); \bar{\partial}, \partial).
$$
There exist two spectral sequences:
\begin{enumerate}
\item $\{E_{r},d_{r}\}$ converges to the total cohomology
$H^{\bullet}(K^{\bullet}(\CH))$ such that
\begin{equation*}
E^{p,q}_{1}=H^{p,q}_{\bar{\partial}}(K^{\bullet}(\CH))
=
H^{q}(K^{p,\bullet}(\CH))
\cong
H^{q}(X,\K_{X,Z}^{p}(\CH));
\end{equation*}
\item $\{\tilde{E}_{r},\tilde{d}_{r}\}$ converges to the total cohomology
$H^{\bullet}(K^{\bullet}(\widetilde{\CH}))$ such that
\begin{equation*}
E^{p,q}_{1}=H_{\bar{\partial}}^{p,q}(K^{\bullet}(\widetilde{\CH}))
=
H^{q}(K^{p,\bullet}(\widetilde{\CH}))
\cong
H^{q}(\widetilde{X},\K_{\widetilde{X},E}^{p}(\widetilde{\CH})).
\end{equation*}
\end{enumerate}
Note that the morphism of double complexes
$
\pi^{\ast}:
K^{\bullet,\bullet}(\CH)
\longrightarrow
K^{\bullet,\bullet}(\widetilde{\CH})
$
induces a morphism of spectral sequences
$
\pi^{\ast}_{r}:E_{r}\rightarrow\tilde{E}_{r}
$
for any $r\geq1$.
Moreover, we have the following result.

\begin{claim}[cf. {\cite[Lemma 4.5]{RYY18}}]
For any $0\leq p, q\leq n$, the pullback of differential forms induces an isomorphism
\begin{equation*}
\pi^{\ast}:
H^{q}(X, \K_{X,Z}^{p}(\CH))
\longrightarrow
H^{q}(\widetilde{X}, \K_{\widetilde{X},E}^{p}(\widetilde{\CH})).
\end{equation*}
\end{claim}

In fact, by \cite[Lemma 4.4]{RYY18} we have $
\pi^{\ast}: \K_{X,Z}^{p}(\CH)
\stackrel{\simeq}\longrightarrow
 \pi_{\ast}\K_{\widetilde{X},E}^{p}(\widetilde{\CH})
$
and  the higher direct image $R^{i}\pi_{\ast} \K_{\widetilde{X},E}^{p}(\widetilde{\CH})=0$ for $i\geq 1$.
Since $\K_{\widetilde{X},E}^{p,\bullet}(\widetilde{\CH})$ is a fine resolution of
$\K_{\widetilde{X},E}^{p}(\widetilde{\CH})$,
the vanishing of the higher direct images implies that
$\pi_{\ast} \K_{\widetilde{X},E}^{p,\bullet}(\widetilde{\CH})$
is a fine resolution of $\pi_{\ast}\K_{\widetilde{X},E}^{p}(\widetilde{\CH})$.
As a result,
we obtain the canonical isomorphisms
\begin{eqnarray*}
\pi^{\ast}: H^{l}(X, \K_{X,Z}^{p}(\CH))
&\stackrel{\simeq}\longrightarrow&
H^{l}(X, \pi_{\ast}\K_{\widetilde{X},E}^{p}(\widetilde{\CH}))\\
&=&
H^{l}(\Gamma(X, \pi_{\ast}\K_{\widetilde{X},E}^{p,\bullet}(\widetilde{\CH})) \\
&=&
H^{l}(\Gamma(\widetilde{X}, \K_{\widetilde{X},E}^{p,\bullet}(\widetilde{\CH})) \\
&=&
H^{l}(\tilde{X}, \K_{\widetilde{X},E}^{p}(\widetilde{\CH})).
\end{eqnarray*}
This means that
$\pi^{\ast}_{1}:E_{1}\rightarrow\tilde{E}_{1}$
is isomorphic, and hence, $\pi^{\ast}_{r}$ is isomorphic for any $r>1$ by a standard result in the spectral sequence theory.
Consequently, the induced morphism \eqref{morph-total} is isomorphic.
This concludes the proof.

\subsection{Proof of Corollary \ref{cor1}}
Recall the Hodge-de Rham spectral sequence
$$
E_{1}^{p,q}=H^{q}(X,\Omega_{X}^{p}(\CH))
\Longrightarrow
H^{p+q}(X, H).
$$
Hence, one has an inequality of the dimension of cohomologies
$$
h^{l}(X,  H)
\leq
\sum_{p+q=l} h^{q}(X,\Omega_{X}^{p}(\CH))
$$
for any $l$.
Then, the Hodge-de Rham spectral sequence degenerates at $E_{1}$
if and only if
$$
h^{l}(X,  H)
=
\sum_{p+q=l} h^{q}(X,\Omega_{X}^{p}(\CH))
$$
for any $0\leq l\leq 2n$.
In fact, by the bundle-valued Dolbeault blow-up formula \cite[Theorem 1.2]{RYY18},
one has
$$
\sum_{p+q=l} h^{q}(\widetilde{X},\Omega_{\widetilde{X}}^{p}(\widetilde{\CH}))
=
\sum_{p+q=l}\big( h^{q}(X,\Omega_{X}^{p}(\CH))
+\sum_{i=1}^{c-1} h^{q-i}(Z,\Omega_{Z}^{p-i}(\iota^{\ast}\CH))\big).
$$
Moreover, according to Theorem \ref{mainresult},
we obtain the following equation
\begin{eqnarray*}
h^{l}(\widetilde{X}, \pi^{-1}H)-
\sum_{p+q=l} h^{q}(\widetilde{X},\Omega_{\widetilde{X}}^{p}(\widetilde{\CH}))
&=&
\big(h^{l}(X, H)-
\sum_{p+q=l} h^{q}(X,\Omega_{X}^{p}(\CH)) \big) \\
&+&
\sum_{i=1}^{c-1}\big(h^{l-2i}(Z, \iota^{-1}H)-
\sum_{p+q=l} h^{q-i}(Z,\Omega_{Z}^{p-i}(\iota^{\ast}\CH)) \big).
\end{eqnarray*}

If the Hodge-de Rham spectral sequence degenerates at $E_{1}$ for $(X, H)$ and $(Z, \iota^{-1}H)$, it follows from the above equation
$$
h^{l}(\widetilde{X}, \pi^{-1}H)=
\sum_{p+q=l} h^{q}(\widetilde{X},\Omega_{\widetilde{X}}^{p}(\widetilde{\CH}))
$$
and thus, the Hodge-de Rham spectral sequence degenerates at $E_{1}$ for $(\widetilde{X}, \pi^{-1}H)$.

Conversely,
suppose the Hodge-de Rham spectral sequence degenerates at $E_{1}$ for $(\widetilde{X}, \pi^{-1}H)$.
Then we have
\begin{eqnarray*}
0&=&\underbrace{h^{l}(X, H)-
\sum_{p+q=l} h^{q}(X,\Omega_{X}^{p}(\CH))}_{\leq 0} \\
&+&
\underbrace{\sum_{i=1}^{c-1} h^{l-2i}(Z, \iota^{-1}H)-
\sum_{p+q=l} h^{q-i}(Z,\Omega_{Z}^{p-i}(\iota^{\ast}\CH)}_{\leq 0}
\end{eqnarray*}
which concludes the proof.


\section{Some remarks}\label{finalrem}

\subsection{Algebraic de Rham cohomology}
Note that the de Rham cohomology theory is not only important for complex manifolds but also significant for smooth algebraic varieties.
One may refer to \cite{Gro66} for the algebraic de Rham cohomology.
Let $Y$ be a smooth projective variety over a field $k$ of characteristic $0$.
Recall that the {\it $l$-th algebraic de Rham cohomology} of $Y$ is defined as the hypercohomology
$$
H_{dR}^{l}(Y/k):=\BH^{l}(Y, \Omega_{Y/k}^{\bullet}),
$$
where $\Omega_{Y/k}^{\bullet}$ is the {\it algebraic de Rham complex} of sheaves of regular differential forms.
Similarly, if $\mathcal{E}$ is an algebraic vector bundle with an integrable connection,
then the K\"{a}hler differential operator on $\Omega_{Y/k}^{\bullet}$ naturally extends to
a differential operator on $\Omega_{Y/k}^{\bullet}\otimes \mathcal{E}$.
The {\it $l$-th twisted algebraic de Rham cohomology} is given by the hypercohomology
$$
H_{dR}^{l}(Y/k; \mathcal{E}):=\BH^{l}(X, \Omega_{Y/k}^{\bullet}\otimes_{\CO_{Y}} \mathcal{E}).
$$

If $k=\CN$, by Serre's GAGA principle, an algebraic vector bundle $\mathcal{E}$ with an integrable connection on $Y$ is equivalent to a holomorphic vector bundle $\mathcal{E}_{an}$ with a flat holomorphic connection on the associated complex manifold $Y_{an}$.
Moreover, there is an isomorphism
$$
H^{l}(Y/\CN; \mathcal{E})\cong H^{l}(Y_{an}; \mathcal{E}_{an}).
$$
As a consequence, Theorem \ref{mainresult} and Serre's GAGA principle imply that the blow-up formula holds for the twisted algebraic de Rham cohomology on smooth projective variety over $\CN$.
Furthermore, the Lesfchetz's principle implies that the blow-up formula holds for the twisted algebraic de Rham cohomology on smooth projective variety over a field of characteristic $0$.

\subsection{A further question}
Let $X$ be a compact complex manifold of dimension $n\geq 2$
and $\CH$ be a holomorphic vector bundle or a flat smooth complex vector bundle.
Given any integer $0\leq s\leq t\leq n$,
the {\it truncated twisted holomorphic de Rham complex},
denoted by $\Omega_{X}^{[s \bullet t]}(\CH)$, is defined to be
$$
\xymatrix{
0\ar[r]& \Omega_{X}^{s}(\CH)
\ar[r]& \Omega_{X}^{s+1}(\CH)
\ar[r]&\cdots
\ar[r]& \Omega_{X}^{t}(\CH)
\ar[r]& 0.
}
$$
Its $l$-th hypercohomology $\BH^{l}(X,\Omega_{X}^{[s \bullet t]}(\CH))$ is a finite dimensional complex vector space.
Then, one may ask the following question.

\begin{quest}\label{final-quest}
Is there a blow-up formula for $\BH^{l}(X,\Omega_{X}^{[s \bullet t]}(\CH))$?
\end{quest}

It seems that this question is quite natural,
since the following two special cases hold:
\begin{enumerate}
\item if $s=t$, it is the {\it bundle-valued Dolbeault blow-up formula} whose $\CH$ is a holomorphic vector bundle (see  \cite[Theorem 1.2]{RYY18});
\item if $s=0$ and $t=n$,  it is the {\it twisted de Rham blow-up formula} whose $\CH$ is a flat $C^{\infty}$ complex vector bundle (see Theorem \ref{mainresult}).
\end{enumerate}

\begin{rem}
To deal with Question \ref{final-quest},
one possible way is to use the same idea as the proof of Theorem \ref{mainresult} and it will be sufficient to obtain the analogous results that of Lemmas \ref{injective-lem}, \ref{bundle-formula} and \ref{key-lem}.
\end{rem}

\section*{Acknowledgement}
The authors would like to thank the Departments of Mathematics of Pennsylvania State University and Universit\`{a} degli Studi di Milano for the hospitalities during their respective visits,
and thank Sheng Rao and Xiangdong Yang for the useful discussions.
In particular, the authors would like to thank the referee for introducing Example \ref{ex-Iwasawa} to them.
This work is partially supported by the NSFC (Grant Nos. 11571242, 11701414), the Science and Technology Research Program
of Chongqing Municipal Education Commission (Grant No.KJ1709216) and the China Scholarship Council.

\end{document}